\documentclass[a4paper,12pt]{article}

\usepackage{amsmath}
\usepackage{amsthm}
\usepackage{amsfonts}
\usepackage{amssymb}

\parskip\smallskipamount
\numberwithin{equation}{section} 

\newtheorem{theorem}{Theorem}

\newtheorem{thm}{Theorem}[section]
\newtheorem{lem}[thm]{Lemma}
\newtheorem{cor}[thm]{Corollary}
\newtheorem{prop}[thm]{Proposition}

\theoremstyle{definition}
\newtheorem{definition}[theorem]{Definition}%chapter

\theoremstyle{remark}
\newtheorem{remark}[theorem]{Remark}

\newcommand{\pr}{\mathbf{P}}
\newcommand{\E}{\mathbf{E}}

\newcommand{\Z}{\mathbb{Z}}

\newcommand{\I}{\mathcal{I}}

\title{Rates of convergence for the three state\\ contact process in one dimension\vspace{-2ex}} 

\author{A. Tzioufas} 
\begin{document}
\maketitle
\vspace{-5ex}

% Your postal address goes here. Actuarial Mathematics and Statistics, Edinburgh EH14 4AS UK,
% E-mail address:\begin{abstract}

%The contact process with infection parameter $\mu$ altered so that first infections of initially susceptible sites occur at rate proportional to $\lambda$ instead is considered. It is known that this process started from a single site infected dies out for $\mu$ less than the contact process' critical value, and survives for $\mu$ greater than that value.
%at rate proportional to
%the time to die out is shown to be exponentially bounded 
%of eventually infected sites 
%the cardinality of

%The basic contact process with infection parameter $\mu$ altered so that infections of never infected sites occur at rate proportional to $\lambda$ instead is considered. It is known that in dimension one the epidemic started from one infected cannot survive when $\mu$ is less than the contact process' critical value, while survival is possible when $\mu$ is greater than that value. In the former case the span of the epidemic is shown to decay exponentially in space and in time. In the latter case and for $\lambda$ less than $\mu$, the ratio of the endmost infected site's velocity to that of the contact process is shown to be no greater than $\lambda / \mu$. 
 %
 %standard
\begin{abstract}
\vspace{-1ex}
The basic contact process with parameter $\mu$ altered so that infections of sites that have not been previously infected occur at rate proportional to $\lambda$ instead is considered. 
Emergence of an infinite epidemic starting out from a single infected site is not possible for $\mu$ less than the contact process' critical value, whereas it is possible for $\mu$ greater than that value. In the former case the space and time infected regions are shown to decay exponentially; in the latter case and for $\lambda$ greater than $\mu$, the ratio of the endmost infected site's velocity to that of the contact process is shown to be at most $\lambda / \mu$.
\end{abstract}

\textit{Key words: Contact processes; immunization }
%(Primary)epidemiological ;

%$\displaystyle{ \alpha_{(\lambda,\mu)} \leq \frac{\lambda}{\mu} \alpha_{(\mu,\mu)}}$, where $\alpha_{(\lambda,\mu)}>0$ is 

%Additionally, 
 % insert the primary Maths Subject Classification number in the first bracket
         % and the secondary ams number(s) in the second bracket
         % e.g. \ams{60E20}{49G03;49F10}
         
         %In this paper the following epidemic model is considered. 
         %= \{\zeta(z), z\in\Z\}
         %Hence, for a configuration $\zeta$, site $x \in \Z$ is in state $\zeta(x)$ .
         %.
\section{Introduction and main results}\label{S0}

%general information about which may be found in 

The \textit{three state contact process} is a continuous-time Markov process $\zeta_{t}$ on the space of configurations $\{-1,0,1\}^{\Z}$ with transition rates corresponding to the following local prescription. Flips of $\zeta_{t}(x)$, the state of site $x \in \Z$ at time $t$,  occur according to the rules: $1 \rightarrow 0$ at rate 1, $-1 \rightarrow 1$  at rate $\lambda \hspace{0.2mm} n(x)$, and $0 \rightarrow 1$  at rate $\mu \hspace{0.2mm}n(x)$, where $n(x)$ takes values $0,1,$ or $2,$ with regard to the number of $y = x-1,x+1$ such that $\zeta_{t}(y)= 1$, and the parameters $\lambda$ and $\mu$ are finite, non-negative constants. The process is an interacting particle system in the common through the literature sense, see for instance \cite{L, D95, L99} and the references therein, although, owing to the inclusion of the third state,  it is not a spin system as the contact process itself. 

%

%originates
The \textit{contact process} and the \textit{standard spatial epidemic} correspond to the particular cases of the three state contact process with $\lambda = \mu$ and $\mu=0$ respectively.
Owing to their simple and elegant definition, these processes serve as basic models for the description of various phenomena and in particular, as their name suggests, for the spread of infections on spatially structured populations. Their introduction in the mathematical literature, along with some fundamental results, originates from Harris in 1974, \cite{H74}, and from the correspondence between Mollison and Kelly in 1977, \cite{M}, respectively. Since then they have been extensively studied and an account of recent developments may be found in \cite{L99} for the former, and in \cite{D88} for the latter, see also \cite{A}.

%These processes have been extensively studied in the mathematical literature since
%both processes possess and thus 

%have been extensively studied in the literature since their introduction from
%.

The three state contact process has been formerly studied by Durrett and Schinazi \cite{DS}, Stacey \cite{S}, and Tzioufas \cite{T,Tzi}, while, independently of \cite{DS}, it was earlier on considered in the physics literature by Grassberger, Chate and Rousseau \cite{GCR}.  The epidemiological interpretation of the process derives from regarding sites in state $1$ as infected, sites in state $0$ as susceptible and previously infected, and sites in state $-1$ as susceptible and not previously infected. This perspective of the process as the generalization of that of the contact process featuring a different initial infection rate has commonly motivated studies in the mathematical literature. The process can thus be thought of as a model for the spread of a disease the primary exposure to which results in permanent alternation of subsequent susceptibility, viz., in the case $\mu < \lambda$ the model confers partial immunization, in contrast to the $\mu > \lambda$ case in which it confers reverse immunization. Theoretically many diseases exhibit the former characteristic that can also be a consequence of imperfect inoculation, while tuberculosis and bronchitis are specific examples of diseases which exhibit the latter one.

%Recently Foss and Zachary \cite{FZ} have provided with the observation that there is some strain of the process in the reformulation of which the corresponding method of proof in \cite{T} can be reapplied.

%, some of  main results originate 

%the first rigorous on these processes 

%Their introduction in the mathematics  while 

%, the survival properties of these processes 

Let $\zeta_{t}^{\eta_{0}}$ denote the process with initial configuration $\eta_{0}$ such that the origin is infected
and all other sites are susceptible and not previously infected. 
The process is said to \textit{survive}
if $\pr(\zeta_{t}^{\eta_{0}} \mbox{ survives}) > 0$, where $\{\zeta_{t}^{\eta_{0}} \mbox{ survives}\}$ is a shorthand for $\{\forall t\geq0, \zeta_{t}^{\eta_{0}}(x) =
1 \mbox{ for some } x \}$, while otherwise it is said to die out. Supposing the process survives and letting $r_{t}$ denote the rightmost infected site in $\zeta_{t}^{\eta_{0}}$, the limit of $t^{-1}r_{t}$ as $t\rightarrow\infty$ on $\{\zeta^{\eta_{0}}_{t} \textup{ survives}\}$ is referred to as the \textit{asymptotic velocity of the rightmost infected}. 
A basic fact about the contact process needed to state our results is the existence of a positive and finite value $\mu_{c}$ at which the following dichotomy occurs. The contact process dies out when the parameter takes values less than the so-called critical value $\mu_{c}$, and survives for values greater than that. (That the process with parameter equal to $\mu_{c}$ dies out, and that this indeed analogously holds on $d$-dimensional lattices for any $d$, although far from straightforward to prove, is well known \cite{BG}).

In the partial immunization case and for any $\mu$ that is less than $\mu_{c}$, as shown in \cite{DS}, the process dies out. Taking a different approach here permits obtaining the following improvement, the first part of which was proved independently of \cite{DS}.  

%This result is improved by the following one by means of 
%this result.
 %Taking a different
%approach  permits showing that the following stronger statement holds.
\begin{thm}\label{thmsub}
For all $\lambda$ and $\mu$ such that $\mu<\mu_{c}$ there exists $\delta<1$ such that
$\pr\left(\exists \hspace{0.5mm} t \textup{ s.t., } \zeta_{t}^{\eta_{0}}(n)=1 \textup{ or } \zeta_{t}^{\eta_{0}}(-n)=1 \right) \leq \delta^{n}, \mbox{ for all } n\geq1$; further, there exist $C$ and $\gamma>0$ such that
$\pr(\exists \hspace{0.5mm} x \textup{ s.t., } \zeta_{t}^{\eta_{0}}(x) =1) \leq C e^{-\gamma t}, \mbox{ for all } t\geq0.$
\end{thm}
%\textup{ for some } x \right
%Remark on generalizations? 
%analogous to Theorem \ref{moninit}
% which is a simple consequence of that the time to die out for this process is known to be exponentially bounded.
%subtlety in the latter arguments is necessary 

%a rigorous formulation of which that suffices for our purposes can also be found below. 
%(A definition of $\beta(\mu)$ that suffices in our case is given below, called the  is a strictly positive constant. 
The following observations in regard to the method of proof of Theorem \ref{thmsub} are in due course here. That of the first part relies on establishing that for subcritical one-dimensional contact processes the probability that the span of infected sites never expands is bounded away from zero uniformly over all finite initial configurations. Ad-hoc arguments which necessitate the uniformity in the first part are used in the proof of the second part for circumventing difficulties stemming from the lack of monotonicity properties of the process. 
%extension for the three state contact process
%(i.e., the interval having the rightmost and leftmost infected sites for endpoints)
%(i.e.\ remains a subset of the initial span for all times)
%due to the lack of monotonicity of the process when $\lambda>\mu$.  %(see Remark \ref{}).instead of known methods (namely, Proposition \ref{geomb} below)

%slightly

Turning to the reverse immunization case and for any $\mu$ that is greater than $\mu_{c}$, 
as shown in \cite{DS} and independently in \cite{TZ}, the process with parameters as such survives. (In fact, Theorem 3 in \cite{DS} establishes the analogous result for the process on the $d$-dimensional lattice for a notion stronger than that of survival). Furthermore, as shown in \cite{T}, the asymptotic velocity of the rightmost infected exists and is a positive constant almost surely. The following comparison result is proved here.
%; note also that when one assumes that $\mu\leq \lambda$ instead this result follows since the process dominates a contact process with parameter $\mu$, as  
%by comparison of the process with the contact process is an immediate consequence of 

% implies in particular that the process for $(\lambda,\mu)$ such  survives . It is a well known result that 
%see e.g.\ )

%the contact process plays a key role in its analysis (see \cite{D80}), and is well known that its is strictly positive when $\mu$ is supercritical. (In fact, it turns out that $\beta(\mu)$ characterizes the critical value, see \cite{D84}, and also \cite{G}). 

\begin{thm}\label{thmsup}
Let $\lambda$ and $\mu$ be such that $\mu>\mu_{c}$ and $\mu > \lambda$. Let  $\alpha$ be the asymptotic velocity of the rightmost infected of the process with parameters $(\lambda,\mu)$, and $\beta$ be that of the contact process with parameter $\mu$. Then, for any such $\lambda$ and $\mu$, $\displaystyle{\alpha \leq (\lambda/ \mu) \beta}$. 
\end{thm}
The proof of this theorem relies on a comparison of the growth of the rightmost infected sites of a sequence of contact processes defined iteratively on the trajectory of the rightmost infected site of the process. Further, known upper bounds of $\beta$ in terms of $\mu$, see p.289 in \cite{L}, yield corresponding upper bounds of $\alpha$ in terms of $\lambda$ and $\mu$ as 
immediate corollaries of this result.

%The  relies on an  

 %, where  of which is based upon the basic result in Lemma \ref{thmpresup}. Remark on the condition $\mu \geq \lambda$ the upper bound given above when combined  provides 

In the following preparatory section the graphical representation is explained and some known results that we use are stated. The remainder of the paper is then devoted to proofs; that of Theorem 1.1 is in Section 3, while that of Theorem 1.2 is given in Section 4.

%We also give some extensions of standard results for contact processes to the three state contact process. %based on known techniques which are extended in this context.  

%The proof is an extention 

 %we believe that an, albeit elementary, proof of this result is of interest.
%; under the assumption that $\mu>c$ we shall refer to the process as $\textit{supercritical}$. 

%In the next section we explain the graphical construction and state some background results; the remainder of the paper is devoted to proofs. and Theorem \ref{thmsup} in 
%Section \ref{Smon} includes the monotonicity results we referred to. 

%is devoted to 
%statements are given in . %\ref{S1} is devoted to  statements. 
%and hence give a proof by use of arguments analogous to those in 

%\section{Converge rates results}

\section{Preliminaries}\label{prel} 
%We introduce..
%\subsection{The graphical construction and background results}\label{grrep} 
 %we will repeatedly use it throughout this paper.
%Let $\lambda$ and $\mu$ be  of the parameters.   
%and are hence ommitted, although  to avoid the possibility of confusion, we note that for proceeding similarly in the $\lambda>\mu$ case, we would).We make abundant use of the graphical representation, indeed in 
%
Graphical representations, also termed as percolation substructures \cite{Griff}, were introduced by Harris \cite{H78} in 1978 and are an important tool in the study of particle systems that aids visualizing their imbeding in space-time by a random graph. Here it is intended for constructing three state contact processes with parameters $(\lambda,\mu)$ as well as contact processes with parameter $\mu$ started from different configurations at different times on the same probability space. Abundant use of this representation is made in the proofs below. %The definition of the process, unlike that of the contact process, is done via an algorithm on the graphical representation and differentiates from that in Section 2 in \cite{S} for brevity and clarity in exposition.  % raphical representation

Suppose that $\mu > \lambda$ and, for all integer $x$ and $y=x-1, x+1$, let $\{T_{n}^{(x,y)}, n\geq1 \}$ and $\{U_{n}^{(x,y)}, n\geq1\}$ be the event times of Poisson processes at rates $\lambda$ and $\mu-\lambda$ respectively. (The case that $\lambda>\mu$ is similar by considering Poisson processes at rates $\mu$ and $\lambda-\mu$ instead). Let also $\{S_{n}^{x}, n\geq 1\}$ be the event times of a Poisson process at rate $1$. All Poisson processes introduced are independent. 

Start with the space-time diagram $\Z \times [0,\infty)$,  
where $\times$ denotes Cartesian product, thought of as giving a time line to each site in $\Z$. The \textit{graphical representation} for parameters $(\lambda,\mu)$ and $\zeta_{t}^{[\eta,s]}$,  $t\geq s$, the \textit{three state contact process} started from a configuration $\eta$ at time $s\geq0$, corresponding to a given realization of the before-mentioned ensemble of Poisson processes are defined as follows. At all times $T_{n}^{(x,y)} = t$ place a directed $\lambda$-\textit{arrow} from $(x, t)$ to $(y, t)$, and, for $t\geq s$, if $\zeta_{t-}^{[\eta,s]}(x)=1$ and $\zeta_{t-}^{[\eta,s]}(y)\in \{-1,0\}$ then set $\zeta_{t}^{[\eta,s]}(y) = 1$. (Here $\zeta_{t-}(x)$ denotes the limit of $\zeta_{t-\epsilon}(x)$ as $\epsilon\rightarrow 0$). At all times $U_{n}^{(x,y)}=t$ place a directed $(\mu-\lambda)$-\textit{arrow} from $(x, t)$ to $(y, t)$, and, for $t\geq s$, if $\zeta_{t-}^{[\eta,s]}(x)=1$ and $\zeta_{t-}^{[\eta,s]}(y)=0$ then set $\zeta_{t}^{[\eta,s]}(y) = 1$. Finally, at all times $S_{n}^{x}=t$ place a \textit{recovery mark} 
at  $(x, t)$, and, for $t\geq s$, if $\zeta_{t-}^{[\eta,s]}(x)=1$ then set $\zeta_{t}^{[\eta,s]}(x) = 0$. Further, $\xi_{t}^{(A, s)}$, $t\geq s$, the (set-valued) \textit{contact process} with parameter $\mu$ started from $A$ at time $s\geq0$, is defined via paths of the graphical representation, which we firstly define as follows. The existence of a connected oriented path from $(x, s)$ to $(y, t)$,  $t\geq s$, that moves along arrows (of either type) in the direction of the arrow and along vertical segments of time-axes without passing through a recovery mark is denoted as $(x, s) \rightarrow (y, t)$, while, that $(x, s) \rightarrow (y, t)$ for some $x \in A$ and $y \in B$, is denoted as $(A, s) \rightarrow (B, t)$. It is then immediate that letting $\xi_{t}^{(A, s)} = \{x: (A, s) \rightarrow (x, t)\}$ gives the correct transition rates (where the equivalence with the configuration-valued contact process can easily be seen by noting that the two types of susceptibility merge when $\lambda=\mu$ and by regarding sites of $\xi_{t}^{(A, s)}$ as infected and others as susceptible). To simplify notation we will write $\zeta_{t}^{\eta}$ for $\zeta_{t}^{[\eta,0]}$ and $\xi_{t}^{A}$ for $\xi_{t}^{(A, 0)}$ and, further, for integer $I$, we write $\xi_{t}^{A} \cap I \not= \emptyset$ instead of $\xi_{t}^{A} \cap \{I\} \not= \emptyset$. 
%$(I, s)$ instead of $(\{I\}, s)$,

%is treated similarly and one can proceed with the definitions above can be carried out analogously. 

% will be coupled analogously and it is then easy to see how to proceeed with the definitions following analogously) 

%and 

A property known as monotonicity in the starting set of the contact process is a well known immediate consequence of the definition by the graphical representation. Here by \textit{monotonicity} we will refer to the following particular form of this property which facilitates its applications below. If a path constrained on $D$ from $(A, s)$ to $(B, t)$ exists, then a path constrained on $D'$ from $(A', s)$ to $(B, t)$ exists for any $D' \supseteq D$ and  $A'\supseteq A$, where a path is said to be constrained on a subset of the integers if it includes vertical segments of time axes of sites of that subset only. This property is also referred to as attractiveness within the literature of interacting particles. The other basic property of the contact process we use is known as self duality. To state it let $(\xi_{t}^{A})$ and $(\xi_{t}^{B})$ denote two contact processes with the same parameter started from $A$ and $B$ respectively, then, 
\begin{equation}\label{selfdual}
\pr(\xi_{t}^{A} \cap B \not= \emptyset) = \pr(\xi_{t}^{B} \cap A \not= \emptyset), 
\end{equation}
for all $t\geq0$. The equality above can be seen to hold by considering paths of the graphical representation that move along time axes in decreasing time direction and along arrows in direction opposite to that of the arrow, and noting that the law of these paths is the same as that of the paths going forward in time defined above. See $\cite{D95}$ and $\cite{L, L99}$ for more information on duality.  

%

%Note that we define and subsequently use the contact process viewed as a set valued process.
%The special marks introduced facilitate? the usual definition of the contact process.
% letting $\eta_{0}$ be the standard initia for any $D$ and $E$ such that $D\supseteq C$ and  $E\supseteq A$l configuration, $\zeta_{t}^{\eta_{0}}$ is denoted as $\zeta_{t}^{O}$, 

%It is important to emphasize that the graphical construction defines all processes $\zeta_{t}^{[\eta,s]}$ and all $\xi_{t}^{A\times s}$ simultaneously on the same probability space, 3\ provides a coupling of all these processes.

% we shall denote $(\zeta_{t}^{[\eta,s]})_{t\geq s}$ as $\zeta_{t}^{[\eta,s]}$. 

%and, letting $\eta_{0}$ be the standard initial configuration, $\zeta_{t}^{[\eta_{0},0]}$ is denoted as $\zeta_{t}^{O}$. Additionally, the event $\{\mathcal{I}(\zeta_{t}^{[\eta,s]}) \not= \emptyset  \mbox{ for all } t\geq s\}$ will be abbreviated below as  $\{\zeta_{t}^{[\eta,s]} \mbox{ survives}\}$. 

%We note that we have produced a version of $\zeta_{t}^{\eta}$ via a countable collection of Poisson processes, this provides well-definedness of the process. Indeed, whenever one assumes that $|\mathcal{I}(\eta)|<\infty$, this is a consequence of standard Markov chains results having an almost surely countable state space; otherwise, this is provided by an argument due to Harris \cite{H}, see Theorem 2.1 in Durrett \cite{D95}.

A miscellany of known results that are used in the proofs is collected together in the remainder of this section. First, an observation regarding monotonicity of the three state contact process which can be found within the last section in \cite{S} is presented. Additional information regarding this property, along with a different and independent proof of this one, can be found in the last chapter of \cite{Tzi}.  
%for any two configurations $\zeta$ and $\zeta'$ we write that , that is, we use the componentwise partial order. 

%for a proof by use of a variation of the basic coupling. 
%
%()
\begin{prop}\label{moninit}
Endow the space of configurations with the natural partial ordering, $\zeta \leq \zeta'$ if and only if $\zeta(x) \leq \zeta'(x)$ for all $x$. Consider the graphical representation for $(\lambda,\mu)$ such that $\mu \geq \lambda$.
If $\eta$ and $\eta'$ are such that $\eta \leq \eta'$, then $\zeta_{t}^{\eta} \leq \zeta_{t}^{\eta'}$, for all $t$. 
\end{prop}

Two well known results for the contact process are presented next; for proofs see \cite{D88,L}. To state them, let $\xi_{t}^{A}$ denote the contact process with parameter $\mu$ started from $A$.

\begin{lem}\label{thmpresup}
Let $R_{t}^{A} = \sup \xi_{t}^{A}$. For any infinite $B$ such that $B \subseteq (-\infty,0]$, $\displaystyle{ \E(R_{t}^{B\cup\{1\}} - R_{t}^{B}) \geq 1}$, for all $t\geq0$.
\end{lem}

For the final statement, recall that $\mu_{c}$ denotes the critical value of the contact process and that $|B|$ denotes the cardinality of a set $B$.  
\begin{thm}\label{thmpre1}
For all $\mu < \mu_{c}$ there exists $\psi>0$ independent of $A$ such that $\pr(\xi_{t}^{A} \not= \emptyset) \leq |A|  e^{- \psi t},$ for all $t \geq 0$. 

%where $|A|$ denotes the cardinality of $A$.
\end{thm}

\section{Proof of Theorem \ref{thmsub}}\label{Ssub} 
The theorem is obtained as a compound of two separate Propositions in this section. To follow the thread of the proof the corresponding remarks in the introduction are useful and, in particular, the reader should bear in mind and notice that the use of uniformity is indispensable in our arguments throughout here.

%  must be kept in mind 
%for digesting the essence of the ideas devised,
%In particular attention must be drawn on
\begin{lem}\label{Tparenthexp}
Let $\hat{\xi}_{t}^{A}$ be the contact process with parameter $\mu$ on $\{\min A,\dots,\max A\}$ started from $A$, $|A|<\infty$. For all $\mu<\mu_{c}$ there exist
$C,\gamma>0$ independent of $A$ such that 
\begin{equation*}\label{cpuni}
\pr\left( \exists \hspace{0.5mm} s\geq t \textup{ s.t., } \hat{\xi}_{s}^{A} \cap \min A \not= \emptyset \textup{ or } \hat{\xi}_{s}^{A} \cap \max A \not= \emptyset \right) \leq C e^{-\gamma t}, \hspace{2mm} \text{ for all } t\geq0.
\end{equation*}

\end{lem}

\begin{proof}
%Consider of the contact process with parameter $\mu$ 
%Consider the graphical representation in which there is only one type of arrows positioned according to event times of Poisson processes at rate $\mu$.

By monotonicity and translation invariance it is sufficient to prove that there exist $C,\gamma>0$ independent of $N \geq 0$ such that
\begin{equation}\label{lemcutpunc}
\pr\left( \exists \hspace{0.5mm} s\geq t \textup{ s.t., } \hat{\xi}_{s}^{[0,N]} \cap  0 \not= \emptyset \textup{ or } \hat{\xi}_{s}^{[0,N]} \cap N \not= \emptyset \right) \leq C e^{-\gamma t}, \hspace{2mm} \text{ for all } t\geq0.
\end{equation}

Define $E_{N,t} = \{ \hat{\xi}_{t}^{[0,N]} \cap N \not= \emptyset \textup{ or } \hat{\xi}_{t}^{[0,N]} \cap 0 \not= \emptyset \}$, $t \geq 0$. We first show that there exists $\psi>0$ such that, for any $N$,
\begin{equation}\label{subdual}
\pr(E_{N,t}) \leq 2 e^{-\psi t}, \hspace{2mm} \mbox{ for all } t\geq0. 
\end{equation}
%consider $(\xi^{[0,N]}_{t})$ and $(\xi^{0}_{t})$, the contact processes parameter $\mu$ started from  $[0,N]$ and $\{0\}$ respectively. 
To this end, we have that there exists a $\psi>0$ such that, for any $N$,
\begin{eqnarray}\label{dualandmon}
\pr\big(\xi^{[0,N]}_{t} \cap N  \not= \emptyset\big) &=& \pr\big(\xi_{t}^{0} \cap [-N,0] \not= \emptyset\big) \nonumber\\
&\leq&  \pr\big(\xi_{t}^{0} \cap \Z \not= \emptyset\big)\hspace{1mm} \leq \hspace{1mm} e^{-\psi t} 
\end{eqnarray}
$t\geq0$, where the equality comes from duality, equation (\ref{selfdual}), and translation invariance, while the two inequalities come from monotonicity and Theorem \ref{thmpre1} respectively. Thus, since $\xi^{[0,N]}_{t}$ stochastically dominates $\hat{\xi}_{t}^{[0,N]}$ by monotonicity, (\ref{subdual}) follows from (\ref{dualandmon}) and translation invariance. 

% from 
%From the follows by monotonicity and . %$E_{N,t} \subseteq \{\xi^{[0,N]}_{t} \cap N  \not= emptyset \textup{ or }\xi^{[0,N]}_{t} \cap 0  \not= \emptyset\}$,   .

For every integer $k\geq1$ define the event $D_{N,k}$ to be such that $\omega \in D_{N,k}$ if and only if $\omega \in E_{N,s}$ for some $s \in (k-1 ,k]$. Because the probability of no recovery mark on the time axes of $N$ and $0$ after the first time $s \in (k-1,k]$ such that $\omega \in E_{N,s}$ and before time $k$ is at least $e^{-2}$, gives that
\begin{equation}\label{expminus}
e^{-2}\pr(D_{N,k}) \leq \pr(E_{N,k})
\end{equation}
for all $k\geq1$.
%where we used that the intersection of $D_{N,k}$ with the event of no recovery mark is a subset of $E_{N,k}$. 

Considering the event $\bigcup\limits_{l\geq0}D_{N, l+\lfloor t \rfloor}$, where $\lfloor \cdot\rfloor$ denotes the floor function, Boole's inequality gives that 
\begin{equation*}
\pr\left( \exists \hspace{0.5mm} s\geq t \textup{ s.t., } \hat{\xi}_{s}^{[0,N]} \cap  0 \not= \emptyset \textup{ or } \hat{\xi}_{s}^{[0,N]} \cap N \not= \emptyset \right) \leq  \sum_{l\geq0} \pr(D_{N, l+\lfloor t \rfloor})
\end{equation*}
$t\geq0$. The proof is thus completed since by (\ref{expminus}) and then (\ref{subdual}) the last display implies (\ref{lemcutpunc}).
\end{proof}

The preceding lemma is used in the proof of the next one as well as in that of Lemma \ref{Teta} below. 
%which in turn is used for the proof of Proposition \ref{propsub2}. 
%We note that the next lemma is used in the proof of Proposition \ref{geomb} following while

\begin{lem}\label{ANScp}
Let $\tilde{\xi}_{t}^{A}$ be the contact process with parameter $\mu$ on $\{\min A-1,\dots,\max A+1\}$ started from $A, |A| < \infty$. For all $\mu<\mu_{c}$ there exists $\epsilon>0$ independent of $A$ such that $\pr\big( \forall \hspace{0.5mm} t\geq0, \tilde{\xi}_{t}^{A} \subseteq [\min A,\max A]\big) \geq \epsilon$.
\end{lem}
%()

\begin{proof}
%Consider the graphical representation in which there is only one type of arrows positioned according to event times of Poisson processes at rate $\mu$.Let $N$ be a non-negative, finite integer.

By monotonicity and translation invariance it is sufficient to show that there exists $\epsilon>0$ independent of $N$ such that
\begin{equation}\label{epsone}
\pr\big( \forall \hspace{0.5mm} t\geq0, \mbox{ }\tilde{\xi}^{[0,N]}_{t} \subseteq [0,N]\big) \geq \epsilon.
\end{equation}

Define $\tilde{E}_{N,t} = \{\tilde{\xi}_{t}^{[0,N]} \cap N+1 \not= \emptyset \textup{ or } \tilde{\xi}_{t}^{[0,N]} \cap -1 \not= \emptyset \}$. We have that there exists $\psi>0$ such that, for any $N\geq0$,
\begin{equation}\label{subdual2}
\pr(\tilde{E}_{N,t}) \leq 2 e^{-\psi t}, \hspace{2mm} \mbox{ for all } t\geq0,
\end{equation}
where (\ref{subdual2}) follows from (\ref{subdual}) by noting that $\tilde{\xi}_{t}^{[0,N]}$ is stochastically smaller than $\hat{\xi}_{t}^{[-1,N+1]}$  by monotonicity, and thus 
$\tilde{E}_{N,t}$ is bounded above in distribution by $E_{N+2,t}$ from translation invariance. (Alternatively, (\ref{subdual2}) can be proved by arguments akin to those used for showing (\ref{subdual}) in the proof of the previous statement). 

%the latter processes association to $\hat{\xi}_{t}^{[0,N+2]}$
%^by letting $\hat{\xi}_{t}^{[0,N]}$ be as in Lemma \ref{Tparenthexp}.  and hence $\pr(\tilde{E}_{N,t}) \leq \pr(E_{N+2,t})$, by .

Define $\tilde{D}_{N,k} = \{ \omega: \omega \in \tilde{E}_{N,s} \mbox{ for some } s \in (k-1,k]\}$, for integer $k\geq1$. Clearly, $\textstyle \bigcap\limits_{k\geq1}\tilde{D}^{c}_{N,k}$ is equal to $\left\{ \forall \hspace{0.5mm} t\geq0, \tilde{\xi}^{[0,N]}_{t} \subseteq [0,N] \right\}$ and  $\pr(\textstyle \bigcap\limits_{k\geq1}\tilde{D}^{c}_{N,k}) = \textstyle \lim\limits_{K \rightarrow \infty} \pr( \bigcap\limits_{k\geq1}^{K}\tilde{D}^{c}_{N,k})$. Thus, Harris' version of the FKG inequality, since the events $\tilde{D}^{c}_{N,1}, \dots, \tilde{D}^{c}_{N,K} $ are all decreasing (see \cite{D88, GRI}), gives that for any $N\geq0$, 
\[
\pr\big(\forall \hspace{0.5mm} t\geq0, \mbox{ }\tilde{\xi}^{[0,N]}_{t} \subseteq [0,N] \big)  \geq  \prod_{k \geq1} \pr(\tilde{D}^{c}_{N,k}).
\] 

However, from $(\ref{subdual2})$ and elementary properties of infinite products we have that there exists $\epsilon>0$ independent of $N$ such that $\prod\limits_{k \geq 1} \big(1 - e\pr(\tilde{E}_{N,k})\big)>\epsilon$. Since also we have that $\pr(\tilde{D}_{N,k}) \leq e^{2} \pr(\tilde{E}_{N,k})$, shown similarly to (\ref{expminus}), the proof is complete from (\ref{epsone}) which thus follows from the last display. 
\end{proof}
%where the last expression is bounded below by $\epsilon>0$ uniformly in $N$ by elementary properties of infinite products, becauseThis completes the proof of . 

%and then (\ref{limsupANk}), such that for any $N\geq0$. 
%from the last display, 

%a fortiori Theorem %\ref{}, (\ref{}).
 %Recall also that we denote by $\mu_{c}$ the critical value of the nearest neighbours one-dimensional contact process. 

We return to consideration of the three state contact process. 
%Occasionally,  we will use the following definition in the rest of this chapter. 

\begin{definition}\label{calI}
Let $\mathcal{I}(\zeta)$ denote the set of infected sites in a configuration $\zeta$, that is,  $\mathcal{I}(\zeta) = \{y \in \Z:\zeta(y) =1\}$. We need to make use of the uniformity over $N$ below in the proof of theorem below. 

\end{definition}

Let $\eta_{N}$ be such that $\mathcal{I}(\eta_{N}) = \{-N,\dots, N\}$ and $\eta_{N}(x)=-1$ for all $x \not\in \mathcal{I}(\eta_{N})$,  $N\geq0$. For $N=0$ the next result reduces to the first part of Theorem \ref{thmsub}. %and note that $\eta_{0}$ is the standard initial configuration.  %that is, $\zeta_{t}^{\eta_{0}}\equiv \zeta_{t}^{\eta_{0}}$.
%Consider $\zeta_{t}^{\eta_{N}}, N<\infty,$ with parameters $(\lambda, \mu)$. 
\begin{prop}\label{geomb}
For all $\lambda$ and $\mu$ such that $\mu<\mu_{c}$ there exists $\epsilon>0$ independent of $N$ such that 
\begin{equation*}\label{eqgeomb}
\pr\Big(\exists \hspace{0.5mm}t \textup{ s.t.}, \mbox{ }\zeta_{t}^{\eta_{N}}(N+n)=1  \textup{ or } \zeta_{t}^{\eta_{N}}(-N-n)=1 \Big) \leq (1-\epsilon)^{n}, \hspace{2mm} \mbox{ for all } n\geq1. 
\end{equation*}

\end{prop}

\begin{proof}
%Consider the graphical representation for $(\lambda, \mu)$ as in the statement. 
%note: if it assymetric put an extra 0 on the end %assosiated
Let $I_{t}^{N} :=  \I(\zeta_{t}^{\eta_{N}})$. 
We first show that there exists $\epsilon>0$ such that, for any $N \geq 0$,
\begin{equation}\label{targ}
\pr\big( \forall\hspace{0.5mm} t \geq 0, \mbox{ }I_{t}^{N} \subseteq  [-N, N] \big)\geq \epsilon. 
\end{equation}
Define the events $B_{N}= \left\{ \forall\hspace{0.5mm}s \in (0,1), I_{s}^{N}\subseteq [-N,N]\right\} \cap \left\{I_{1}^{N} \subseteq [-N+1,N-1] \right\}$ and $F_{N} = \{\forall\hspace{0.5mm} t \geq 1, \mbox{ }I_{t}^{N} \subseteq  [-N+1, N-1]\}$. Since $\{\forall\hspace{0.5mm} t \geq 0, \mbox{ }I_{t}^{N} \subseteq  [-N, N] \}\supseteq F_{N} \cap B_{N}$ and, by Lemma \ref{ANScp} and the Markov property at time $1$, the $\pr\big(F_{N}| B_{N} \big)$ is uniformly in $N$ bounded away from zero, it is sufficient to show that $\pr(B_{N})$ also is. For this consider the event $B_{N}'$ that: a) for all times $s \in (0,1)$ no arrow exists from $(N, s)$ to $(N+1,s)$ as well as from $(-N, s)$ to $(-N-1, s)$, b) a recovery mark exists on the time axis of $N$ within $(0,1])$ and $-N$ on $(0,1]$, and, c) no arrow exists from $(N-1, s)$ to $(N,s)$ and over $(-N+1, s)$ to $(-N, s)$, for all times $s \in (0,1]$. Note that b) implies that there is a $t\in(0,1]$ such that $I_{t}^{N} \subseteq [-N+1,N-1]$ and c) assures that this holds for $t=1$, and hence by a) we have that $B_{N} \supseteq B_{N}'$. This proves (\ref{targ}) because $B_{N}'$ has strictly positive probability which is independent of $N$ from translation invariance.   

From (\ref{targ}) and monotonicity (of the contact process) we have that indeed for any $\eta$ such that $\eta(x) \not=-1$, $\forall \hspace{0.1mm} x \in [\min\I(\eta), \max\I(\eta)]$, $\pr\big( \forall \hspace{0.5mm}t \geq 0,\mbox{ } \I(\zeta_{t}^{\eta}) \subseteq  \I(\eta) \big)\geq \epsilon$, and the proof is completed by repeated applications of the Strong Markov Property. 
%for after the recovery of $N$ and $-N$, as in (b), up until time 1). 

%From this and  of the process, we have that 
%because by , we have that there exists $\epsilon_{1}>0$ such that for any $N$,
%\begin{equation}\label{epsiSpan}
%\pr\left(\forall\hspace{0.5mm}s \geq 1, I_{s}^{N} \subseteq [-N+1,N-1] \vline \mbox{ }E_{N}\right)\geq \epsilon_{1},
%\end{equation}
%it thus suffices to show that there is an $\epsilon_{2}>0$ such that for any $N$, 
%$\pr(E_{N}) \geq \epsilon_{2}$; however, for any $N$, $E_{N}$ contains an event formulated in terms of the Poisson processes associated with the sites at $N,N-1,-N+1 \mbox{ and} -N$, up to time $1$, thus the proof is complete. 

%To see that this is sufficient, note that by the nearest neighbours assumption and the strong Markov property, letting $B$ be any interval such that $|B| = |A| +n$, by monotonicity of the contact process, we also have that
%\begin{eqnarray*}
%\pr(\tau_{n+1}^{\eta_{A}} < \infty| \mbox{ }\tau_{n}^{\eta_{A}} < \infty) &\leq& \pr(\tau_{1}^{\eta_{B}} < \infty) \\
% &\leq&  (1-\epsilon)
%\end{eqnarray*}
%for all $n\geq0$; iterrating the last display gives (\ref{eqgeomb}). 

%that by time $1$, (a) no arrow exists from $N$ to $N+1$ nor from $-N$ to $-N-1$ respectively, (b) a recovery mark occurs at $N$ and at $-N$ and, (c) no attempt of infection occurs to $N$ and to $-N$ from $N-1$ and from $-N+1$ respectively, 
\end{proof}

In the proof of Proposition \ref{propsub2} below we need to use the preceding proposition as well as the next corollary. To state the latter, let $H$ be the collection of configurations $\eta$ such that $|\I(\eta)|<\infty$ and $\eta(x) \not=-1, \mbox{ } \forall x \in [\min\I(\eta), \max\I(\eta)]$, and further define the stopping time $T^{\eta}:= \inf\{t\geq0 : \I(\zeta_{t}^{\eta}) \not\subseteq [\min\I(\eta), \max\I(\eta)]\}$, $\eta \in H$. Regarding notation, $1_{E}$ denotes the indicator of event $E$ throughout.   
%Consider $\zeta_{t}^{\eta}, \eta \in H$, with parameters $(\lambda, \mu)$.
\begin{lem}\label{Teta}
For all $\lambda$ and $\mu$ such that $\mu<\mu_{c}$ there exist $C$ and $\theta>0$ independent of $\eta \in H$ such that $\E(e^{\theta T^{\eta}1_{\{T^{\eta}<\infty\}}}) \leq C$.
\end{lem}

\begin{proof}
This follows from Lemma \ref{Tparenthexp}  by the integral representation of expectation since, for any $\eta \in H$, $\{t \leq T^{\eta} <\infty\}$ is bounded above in distribution by $\{\exists \hspace{0.5mm} s\geq t \textup{ s.t., } \hat{\xi}_{s}^{\I(\eta)} \cap \min \I(\eta) \not= \emptyset \textup{ or } \hat{\xi}_{s}^{\I(\eta)} \cap \max \I(\eta) \not= \emptyset\}$.
\end{proof}

%The proof of the next statement completes the section 

Consider $\zeta_{t}^{\eta_{0}}$ with parameters $(\lambda,\mu)$ and let $I_{t} = \I(\zeta_{t}^{\eta_{0}})$.
The final statement of this section is the second part of Theorem \ref{thmsub}. %Recall that the process started from the standard configuration is denoted by $\zeta_{t}^{\eta_{0}}$. 

\begin{prop}\label{propsub2}
For all $\lambda$ and $\mu$ such that $\mu<\mu_{c}$ there exist $C$, $\gamma>0$ such that
$\pr(I_{t} \not= \emptyset) \leq C e^{-\gamma t},$ for all $t\geq0$. 
\end{prop}

%Note that $\tau_{K-1}$ is a , 

\begin{proof} 
%Consider the graphical representation for $(\lambda, \mu)$ as in the statement. 
% where by convention $\min \emptyset = \infty$,
%(which is not a stopping time)
%we have that
%, where by convention $\min \emptyset = \infty$.

Let $S_{t} =  [\min I_{t} , \max I_{t} ] \cap \Z$; define the stopping times $\tau_{k} = \inf\{t\geq0: |S_{t}| = k\}$, $k\geq1$; define also $K = \inf\{k: \tau_{k} = \infty\}$, and further $\sigma_{K} = \inf\{s\geq 0 : I_{s+\tau_{K-1}} = \emptyset\}$. Clearly $\{I_{t} \not= \emptyset\}$ equals $\{\tau_{K-1}+\sigma_{K} \geq t\}$, thus, showing that $\tau_{K-1}$ and $\sigma_{K}$ are exponentially bounded implies the statement since the sum of two exponentially bounded random variables is itself exponentially bounded (where, a simple proof of this fact can be done by using in turn the integral representation of expectation, the Chernoff bound and the Cauchy-Schwartz inequality).  Towards this, because $K$ is exponentially bounded by Proposition $\ref{geomb}$ and by set theory we have that, for all $a>0$,
\begin{equation*}\label{tauK}
\pr(\tau_{K-1} > t) \leq \pr(K > \lceil at \rceil ) + \pr(\tau_{K-1} > t, K \leq \lceil at \rceil ), 
\end{equation*}
$t\geq0$, it suffices to show that (i) there is $a>0$ such that $\tau_{K-1}$ is exponentially bounded on $\{K \leq \lceil at \rceil\}$ and, by repeating the argument in the last display, that (ii) $\sigma_{K}$ is exponentially bounded on $\{K \leq \lceil t \rceil\}$. % and similarly for $\sigma_{K}$. %to that in the last display above  
%To see this note that %we will show that:  and,   that led to the inequality of for fixed $a=1$, 

%It is sufficient to prove that, the non-stopping time, $\tau_{K-1}$ and $\sigma_{K}$ are both exponentially bounded. To see this note that then, equivalently, there are $\theta_{1},\theta_{2}>0$ such that $\E(e^{\theta_{1}\tau_{K-1}})$ and $\E(e^{\theta_{2}\sigma_{K}})$ are finite. Thus, because $\{I_{t} \not= \emptyset\}$ is equal to $\{\tau_{K-1}+\sigma_{K} \geq t\}$ we have that, for all $\gamma>0$, 
%\begin{eqnarray*}\label{twocomp}
%\pr(I_{t} \not= \emptyset) &\leq& e^{-\gamma t} \E\left(e^{\gamma (\tau_{K-1}+\sigma_{K})}\right) \\
%& \leq& e^{-\gamma t} \E\big(e^{\frac{\gamma}{2}\tau_{K-1}}\big)\E\big(e^{\frac{\gamma}{2} \sigma_{K}}\big)
%\end{eqnarray*}
%$t\geq0$, where the second inequality is due to the Cauchy--Schwartz inequality, and the proof is complete by choosing $\gamma <2\min\{\theta_{1}, \theta_{2}\}$ in last display .

%, which is a consequence of the nearest neighbours assumption, 
Towards (i), let $H$ and $C,\theta>0$ be as in Lemma \ref{Teta}. By the Strong Markov Property and because $\zeta_{\tau_{k-1}}^{\eta_{0}} \in H$, we have that
\begin{eqnarray}\label{subindu}
\E( e^{\theta \tau_{k}1_{\{\tau_{k} < \infty\}}}) &\leq& \E(e^{\theta \tau_{k-1}1_{\{\tau_{k-1} < \infty\}}}e^{\theta (\tau_{k}-\tau_{k-1})1_{\{(\tau_{k}-\tau_{k-1}) < \infty\}}}) \nonumber \\
&\leq& C \E(e^{\theta \tau_{k-1}1_{\{\tau_{k-1} < \infty\}}})  \nonumber
\end{eqnarray}
$k\geq1$, which by iteration gives that $\E( e^{\theta \tau_{k}1_{\{\tau_{k} < \infty\}}}) \leq C^{k}$. Using this and set theory gives that, for all $a>0$,
\begin{eqnarray*}\label{tauN2}
\pr(\tau_{K-1} > t, K \leq \lceil at \rceil ) &\leq& \sum_{k=1}^{\lceil at \rceil}e^{-\theta t}\E( e^{\theta \tau_{k-1}1_{\{\tau_{k-1} < \infty\}}}) \nonumber\\
&\leq& \lceil at \rceil e^{-\theta t} C^{\lceil at \rceil},
\end{eqnarray*}
$t\geq0$, and the claim follows from the last display by choosing $a>0$ such that $e^{-\theta}C^{\lceil a \rceil}<1$. 

Towards (ii), let $\hat{\xi}_{t}^{[1,k]}$ denote the contact process with parameter $\mu$ on $\{1,\dots,k\}$ started from all sites infected, it then follows from Theorem \ref{thmpre1} that the $\sum_{k=1}^{\lceil t \rceil} \pr(\sigma_{k}>t, K=k )$ is exponentially bounded in $t$, since $\{\sigma_{k}1_{\{K=k\}}\geq t\}$ is stochastically bounded above by $\{\hat{\xi}_{t}^{[1,k]} \not= \emptyset\}$.  

%by monotonicity,
 
%, 
%, and the proof is thus completed. 
\end{proof}

\begin{remark}
\textup{ It follows from Proposition \ref{geomb} and bounded dominated convergence that $\E|\zeta_{t}^{\eta_{0}}| \rightarrow 0$, as $t\rightarrow \infty$. 
Neither the technique of the proof of Theorem 6.1 in \cite{GRI} nor that of Proposition 1.1 in \cite{AJ} adapt to extend this conclusion to Proposition \ref{propsub2} due to lack of properties of $\zeta_{t}$ analogous to monotonicity and (sub)additivity of the contact process respectively.}
\end{remark}

%This combined with a subadditivity property and an argument as that for the contact process in e.g.\ , or an adaptation of the method of the proof of , would imply Proposition \ref{propsub2}. However neither of these approaches seems to apply here, since neither additivity nor (the weaker) monotonicity property of contact processes are known to generalize in this case

%it appears that   due to the lack of a generic monotonicity (in the initial configuration) of additivity property for the process when $\lambda >\mu$ mentioned(?).

%. \eta_{0}vercoming this difficulty  for
%For choosing the approach in the proof instead of  combined with the consequence of  is crucial.

%Because by coupling
\section{Proof of Theorem \ref{thmsup}}\label{Ssup} 

Let $\zeta^{\bar{\eta}}_{t}$ be the three state contact process with parameters $(\lambda, \mu)$ and initial configuration $\bar{\eta}$ such that $\bar{\eta}(x) =1$ for all $x \leq 0$ and $\bar{\eta}(x) =-1$ for all $x\geq1$. Let also $\bar{r}_{t} = \sup\I(\zeta^{\bar{\eta}}_{t})$ and $\bar{x}_{t} = \sup_{s\leq t} \bar{r}_{s}$. In this section we concentrate on the study of $\bar{r}_{t}$, where the necessary connection between $\zeta^{\bar{\eta}}_{t}$ and $\zeta_{t}^{\eta_{0}}$ for establishing Theorem \ref{thmsup} is given by Corollary \ref{asconv} below. The following lemma is required in the latter's proof. %the proof of which is based on an application of the  and goes through adapting the argument in the  corresponding proof regarding the right endpoint of the contact process, see e.g.\ Theorem 2.19 in \cite{L}.

%The proof of the next lemma below 

\begin{lem}\label{subadd}
If $\mu\geq\lambda$, then $\displaystyle{ \frac{\bar{x}_{n}}{n} \rightarrow a}$ almost surely, where  $\displaystyle{  a =  \inf_{n\geq 0} \frac {\E(\bar{x}_{n})} {n}}$ and $a \in [-\infty,\infty)$. If additionally $a> -\infty$, then $\displaystyle{ \frac{\bar{x}_{n}}{n} \rightarrow a}$ in $L^{1}$.% that is, $\displaystyle{ \lim_{n \rightarrow \infty} \vline \frac {\E(\bar{x}_{n})} {n} - a \vline = 0}$.
\end{lem}

\begin{proof}

%to see this note that $\bar{x}_{0,u}$ will be maximum if at time $s$ all vertices $ x \leq \bar{x}_{s} $ are infected and this equals the left side, i.e. $x_{t}$ is subadditive.  all processes defined in the proof are coupled.

%Consider the graphical representation for any $(\lambda,\mu)$ such that $\mu\geq\lambda$. 
%, .  
%the configuration%and consider the coupled process $\zeta_{t}^{[\eta_{\bar{x}_{s}},s]}$
% $\bar{x}_{0} = 0$
Let $\eta_{y}$ denote the configuration such that $\eta_{y}(z) =1$ for all $z\leq y$, and $\eta_{y}(z) =-1$ for all $z\geq y+1$. For any times $s$ and $u$ such that $s \leq u$, define 
\[
\bar{x}_{s,u} = \max\{y: \zeta_{t}^{[\eta_{\bar{x}_{s}},s]}(y)=1, \mbox{ for some } t \in [s,u] \} - \bar{x}_{s}, 
\] 
where, note that, $\bar{x}_{0,u} = \bar{x}_{u}$. We aim to show that $\{\bar{x}_{m,n}, m \leq n\}$ satisfies the conditions of the subadditive ergodic theorem. We have that
\[
\mbox{a) } \hspace{5mm} \bar{x}_{0,s} + \bar{x}_{s,u} \geq \bar{x}_{0,u}
\]
since, by monotonicity in the initial configuration, Proposition \ref{moninit},  $\zeta_{t}^{[\eta_{\bar{x}_{s}},s]} \geq \zeta^{\bar{\eta}}_{t}$, for all $t\geq s$. 
We further have that $\bar{x}_{s,u}$ is equal in distribution to $\bar{x}_{0,u-s}$ and is independent of $\bar{x}_{0,s}$ by translation invariance and independence of Poisson processes at disjoint parts of the graphical representation respectively. Thus, 
\[\mbox{b)} \hspace{5mm} \{\bar{x}_{(n-1)k,nk},n \geq 1\} \mbox{ are i.i.d.\ for each }  k\geq1, \]
%in particular stationary and ergodic, 
and, furthermore,
\[\mbox{c)} \hspace{5mm} \{\bar{x}_{m,m+k} , k \geq 0\} =  \{\bar{x}_{m+1,m+k+1} , k \geq 0\} \mbox{ in distribution,  for each } m\geq1.\]
By ignoring recovery marks in the representation, $\bar{x}_{t}$ is bounded above in distribution by the number of arrivals of a Poisson process at rate $\lambda$ in $(0,t]$, and thus, from standard properties of Poisson processes we also have that
\[
\mbox{d)} \hspace{5mm} \E(\max\{\bar{x}_{0,1},0\} ) < \infty.
\]  
The result's statement follows from the conclusion of Theorem 2.6, Chapter VI in \cite{L}, since 
the conditions under which it holds correspond to a)--d) above.
\end{proof}

Let $\lambda$ and $\mu$ be such that $\mu>\mu_{c}$ and $\mu > \lambda$ and, further, let $\alpha>0$ be the corresponding value of the asymptotic velocity of the rightmost infected.
The next statement is obtained based on results in \cite{T}.

% where, recall from the introduction, that $\mu_{c}$ denotes the contact process' critical value. Further recall that the , , exists and is strictly positive. 
% Let be the asymptotic velocity of the rightmost infected of the process with parameters $(\lambda,\mu)$.
%We then have that 

\begin{cor}\label{asconv}
$\displaystyle{\frac{\E\bar{r}_{t}}{t} \rightarrow \alpha}$. 
\end{cor}

%appealing to
\begin{proof}
From the embedding of processes started from configurations with one infected site and all other susceptible and not previously infected on the trajectory of $\bar{r}_{t}$ as explicitly done in the statement of Lemma 4.4 in \cite{T} we immediately have that 
\begin{equation}\label{irbar}
\frac{\bar{r}_{t}}{t} \rightarrow \alpha, \mbox{ almost surely}, 
\end{equation}
since both $Y_{N}$ and $T_{Y_{N}}$ in that statement are almost surely finite by Proposition 4.2 in the same paper. From Lemma \ref{subadd} and because $\bar{x}_{n} \geq \bar{r}_{n}$, the last display gives that $\displaystyle{ \frac{\bar{x}_{n}}{n} \rightarrow a}$ in $L^{1}$ for $a>0$, so that by the direct part of the theorem in section 13.7 in \cite{W} it follows that $\displaystyle{ \bar{x}_{n} / n}$ are uniformly integrable and thus, using that $\bar{r}_{n} \leq \bar{x}_{n}$ again,  $\displaystyle{ \bar{r}_{n} / n}$ also are. The latter along with (\ref{irbar}) imply from the reverse part of the before-mentioned theorem in \cite{W} that $\displaystyle{ \frac{\bar{r}_{n}}{n}\rightarrow \alpha}$ in $L^{1}$. The extension along real times then comes elementarily by using that $\max\limits_{t\in (n,n+1]}(\bar{r}_{t} -\bar{r}_{n})$ and  $\max\limits_{t\in (n,n+1]}(\bar{r}_{n+1} -\bar{r}_{t})$ are bounded above in distribution by the number of arrivals of a Poisson process at rate $\mu$ in $(0,1]$. 
%$\displaystyle{\frac{\E\bar{r}_{n}}{n} \rightarrow \alpha}$.
% thus obtaining the result for integer times.

%all processes defined in the proof are coupled.
%Finally by the conclusion of the before mentioned theorem, we also have that $(\bar{x}_{t}/t, t\geq0)$ is a collection of uniformly integrable random variables, by appealing to a standard necessary and sufficient condition for this property to hold, see paragraph 13.7 of Williams \cite{W}. Let $\bar{x}_{t} = \sup_{s\leq t} \bar{r}_{s}$. 
%consider the graphical representation for $(\lambda,\mu)$ such that $\mu\geq\lambda$ and $\mu>\mu_{c}$.

%This completes the proof by standard results of Poisson processes.
%and the first Borel-Cantelli lemma exactly as it is done for the contact process, see Theorem 2.19 of Liggett (1985). 
\end{proof}

%this completes the proof. %we have that $\displaystyle{ \lim_{t\rightarrow \infty} \frac{\E\bar{r}_{t}}{t} =  \alpha}$ and also $\displaystyle{\lim_{t\rightarrow \infty} \frac{\E R_{t}}{t} = \beta}$.
%(i.e. $\zeta^{\bar{\eta}}_{t}(\bar{r}_{t}+1)=-1$, i.e.\ $t$ is such that $\bar{r}_{t} = \sup_{s\leq t}\bar{r}_{s}$. 
%In words, we and consider the times, at each of those times consider the set of infected sites of $\zeta_{t}^{\bar{\eta}}$ and start off a coupled contact process with this as starting set and parameter $\mu$. 

%$(\xi_{t}^{n}; n\geq0)$
\begin{proof}[proof of Theorem \ref{thmsup}]
%A sequence of contact processes parameter $\mu$ will be defined on the trajectory of $\bar{r}_{t}$.
%

% proof's setting (or set up) is given with commentary for expository purposes first. %Its iteration immediately afterwards comprises 
Let $\xi_{t}^{0}, t\geq0,$ be the contact process with parameter $\mu$ such that $\xi_{0}^{0}= \{\dots, -1, 0\}$ and let also $R_{t}^{0} = \sup \xi_{t}^{0}$. We prove the following stronger statement 
\begin{equation}\label{strbnd}
\E \bar{r}_{t} \leq \frac{\lambda}{\mu}\mbox{ } \E R_{t}^{0},
\end{equation}
for all $t \geq0$, which implies the result from Corollary \ref{asconv}. 
 
The first step of the iterative definitions following is outlined with remarks for purposes of illustration. By coupling, $\bar{r}_{t} = R_{t}^{0}$ for all $t$ up until the first time $s$ such that $\bar{r}_{s} = \bar{x}_{s}$ and a $(\mu-\lambda)$-arrow exists from $\bar{r}_{s}$ to $\bar{r}_{s}+1$. Observe that the rightmost infected of the contact process started at time $s$ from $\I(\zeta^{\bar{\eta}}_{s})$ coincides with $\bar{r}_{t}$ up until the first time 
$u$, $u > s$,  at which $\bar{r}_{u} = \bar{x}_{u}$ and a $(\mu-\lambda)$-arrow from $\bar{r}_{u}$ to $\bar{r}_{u}+1$ is present, and further observe that $\I(\zeta^{\bar{\eta}}_{s})$, the starting set of this contact process, equals $\xi_{s}^{0} \backslash \sup\xi_{s}^{0}$.

Define iteratively the stopping times
\begin{equation}\label{eq:upsilons}
\upsilon_{n} = \inf\{t \geq \upsilon_{n-1}: R_{t}^{n-1}= \bar{r}_{t} + 1\},
\end{equation}
where $\upsilon_{0}=0$ and $n\geq1$; define further $\xi_{t}^{n}:= \xi_{t}^{(\I(\zeta^{\bar{\eta}}_{\upsilon_{n}}), \upsilon_{n})}$, $t\geq \upsilon_{n},$ and $R_{t}^{n} = \sup \xi_{t}^{n}$. 
Then, 
\begin{equation}\label{Rtn}
 \bar{r}_{t} = R_{t}^{n-1}, \mbox{ for all } t \in [\upsilon_{n-1},\upsilon_{n}), 
\end{equation}
%and also that
\begin{equation}\label{Hupsilsup}
\xi_{\upsilon_{n}}^{n-1} = \xi_{\upsilon_{n}}^{n} \cup \{\bar{r}_{\upsilon_{n}}+1\}, \mbox{ for all } n\geq1,
\end{equation}
which can be seen to hold from the first and second observation respectively in the outline above.
Define also $F_{t} = \sup\{n: \upsilon_{n} \leq t\}$. We will show that
\begin{equation}\label{fracpunch}
\E(F_{t})= \frac{\mu-\lambda}{\lambda} \E(\bar{x}_{t})
\end{equation}
and, further, that 
\begin{equation}\label{RFt}
\E(R_{t}^{0}-\bar{r}_{t}) \geq \E F_{t}
\end{equation}
$t \geq 0$. 
Note that, since $\bar{x}_{t} \geq \bar{r}_{t}$, (\ref{fracpunch}) gives that $\displaystyle{ \E(F_{t}) \geq \frac{\mu-\lambda}{\lambda} \E(\bar{r}_{t})}$,  which, combined with (\ref{RFt}), implies (\ref{strbnd}). Thus, showing the last two displays above gives (\ref{strbnd}) from which the proof is complete.

%Showing completes the proof since then  follows. To see this  .   %it is sufficient to prove that, for all

%note that then $R_{t}=\bar{r}_{t}$ for all $t \leq \upsilon_{1}$, and also that $\bar{r}_{t}$ will then agree with the right endpoint of the contact process with parameter $\mu$ started from $\I(\bar{\zeta}_{\upsilon_{1}}) \times \upsilon_{1}$ from time $\upsilon_{1}$ up until $\upsilon_{2}$. 
%the contact process with parameter $\mu$ started from $\I(\bar{\zeta}_{\upsilon_{n}})$ at time $\upsilon_{n}$. 
 
% 
%further, . In words 

%note that, by elementary properties of the Poisson processes used for the construction  $\upsilon_{n}<\infty$ almost 

%, see (\ref{propermon})
Let $\mathcal{F}_{t}$ denote the sigma algebra associated to the Poisson processes in the graphical representation up to time $t$ and recall that $1_{E}$ denotes the indicator of event $E$. We first prove (\ref{RFt}).  From $(\ref{Rtn})$ we have that $R^{0}_{t}-\bar{r}_{t} = \sum_{n=1}^{\infty}(R^{n-1}_{t}-R_{t}^{n}) 1_{\{F_{t} \geq n\}}$. This and the monotone convergence theorem, which applies because  $R_{t}^{n-1} \geq R_{t}^{n}$ by monotonicity of the contact process, give that 
\begin{equation}\label{ER0minusr}
\E(R_{t}^{0}-\bar{r}_{t}) = \sum_{n=1}^{\infty}\E\left( (R^{n-1}_{t}-R_{t}^{n}) 1_{\{F_{t} \geq n\}} \right),
\end{equation}
$t\geq0$. Further, Lemma \ref{thmpresup} and $(\ref{Hupsilsup})$ by use of the Strong Markov Property give that
 %for $B$ there taken to be $\xi_{\upsilon_{n}}^{n}$
\begin{equation}\label{ER0minusr1}
\E\left( (R^{n-1}_{t}-R_{t}^{n}) 1_{\{F_{t} \geq n\}} \right) \geq \pr(F_{t} \geq n), 
\end{equation}
$n\geq 1$, where we used that ${\{F_{t} \geq n\}} = \{\upsilon_{n} \leq t\} \in \mathcal{F}_{\upsilon_{n}}$. Thus (\ref{RFt}) follows by plugging (\ref{ER0minusr1}) into (\ref{ER0minusr}) and the telescopic formula for expectation. 
%if $M=$ then ;  while if $M= \tilde{S_{1}}$ then $\bar{r}_{\tilde{T_{1}}} = \bar{x}_{\tilde{T_{1}}}=-1$.

Towards (\ref{fracpunch}) some additional definitions are necessary. Recall the setting of the graphical representation from Section \ref{prel}. Let $\tilde{T}_{1}:= T_{1}^{(0,1)}$, $\tilde{S}_{1}:= S_{1}^{0}$, $\tilde{U}_{1}:= U_{1}^{(0,1)}$ and also 
define the events $A_{1} = \{\min\{\tilde{T}_{1}, \tilde{S}_{1}, \tilde{U}_{1}\} = \tilde{U}_{1}\}$ and $B_{1}= \{\min\{\tilde{T}_{1}, \tilde{S}_{1}, \tilde{U}_{1}\}=\tilde{T}_{1}\}$. At time $\tau_{0}:= 0$ the first competition takes place in the sense that on $A_{1}$, $\bar{r}_{\tilde{U}_{1}}= 0$ and $R_{\tilde{U}_{1}}^{0}=1$ (and hence $\upsilon_{1} = \tilde{U}_{1}$); while on $B_{1}$, $\bar{r}_{\tilde{T}_{1}} = \bar{x}_{\tilde{T}_{1}}=1$. We repeat these inductively as follows. For all $n\geq1$ consider
\[
\tau_{n} = \inf\{ t\geq \min\{\tilde{T}_{n}, \tilde{S}_{n}, \tilde{U}_{n}\}: \bar{r}_{t} = \bar{x}_{t}\}, 
\]
and let $\tilde{T}_{n+1}=\inf\limits_{k\geq1}\{T_{k}^{(\bar{r}_{\tau_{n}},\bar{r}_{\tau_{n}}+1)}:T_{k}^{(\bar{r}_{\tau_{n}},\bar{r}_{\tau_{n}}+1)}> \tau_{n}\}$, i.e.\ the first time a $\lambda$-arrow exists from $\bar{r}_{\tau_{n}}$ to $\bar{r}_{\tau_{n}}+1$  after $\tau_{n}$, and $\tilde{U}_{n+1} = \inf\limits_{k\geq1}\{U_{k}^{(\bar{r}_{\tau_{n}},\bar{r}_{\tau_{n}}+1)}: U_{k}^{(\bar{r}_{\tau_{n}},\bar{r}_{\tau_{n}}+1)}> \tau_{n}\}$,  i.e.\ the first such time a $(\mu-\lambda)$-arrow exists, and further $\tilde{S}_{n+1} =\inf\limits_{k\geq1} \{S_{k}^{\bar{r}_{\tau_{n}}}: S_{k}^{\bar{r}_{\tau_{n}}} >  \tau_{n}\}$, i.e.\ the first time that a recovery mark exists on $\bar{r}_{\tau_{n}}$ after $\tau_{n}$. Define also the events $A_{n+1}:=  \{\tilde{U}_{n+1} < \min\{\tilde{T}_{n+1},\tilde{S}_{n+1}\}\}$ and $B_{n+1}:= \{ \tilde{T}_{n+1} < \min\{ \tilde{U}_{n+1}, \tilde{S}_{n+1}\}\}$. The stopping times $\tau_{n}$ can be thought of as the time that the $n+1$ competition, in the sense explained above, takes place.

Letting $N_{t} = \sup\{n: \tau_{n} <t\}$, we have that $\bar{x}_{t} = \sum\limits_{n=1}^{N_{t}}1_{B_{n}}$ and also that $F_{t} = \sum\limits_{n=1}^{N_{t}} 1_{A_{n}}$, where the latter can be seen by noting that $\upsilon_{n}$ can also be expressed as the first $\tilde{U}_{k}$ after $\upsilon_{n-1}$ such that $\tilde{U}_{k} < \min\{ \tilde{T}_{k}, \tilde{S}_{k}\}$. The last two equalities and assuming that $\E(N_{t})<\infty$ imply  (\ref{fracpunch}) as follows. Since conditional on $\zeta^{\bar{\eta}}_{\tau_{n}}$ the events $A_{n+1}$ and $B_{n+1}$ are independent of $\{N_{t} \geq n+1\}= \{N_{t} \leq n\}^{c} \in \mathcal{F}_{\tau_{n}}$ from the Strong Markov Property, emulating the proof of Wald's lemma and then using a basic result about competing Poisson processes gives that $\displaystyle{\E(\bar{x}_{t}) = \E(N_{t})\frac{\lambda}{\mu+1}}$, and also that $\displaystyle{\E(F_{t}) = \E(N_{t}) \frac{\mu-\lambda}{\mu+1}}$, hence, (\ref{fracpunch}) follows by combining these last two equalities.

It remains to show that $\E(N_{t})<\infty$. Ignoring recovery marks gives that $R^{0}_{t}$ is bounded above (in distribution) by $\Lambda_{\mu}[0,t)$,  the number of arrivals of a Poisson process at rate $\mu$ in $[0,t)$, and further that $\bar{x}_{t}$ is bounded above by $\Lambda_{\lambda}[0,t)$, while also $D_{t}$, the total number of recovery marks  on the trajectory of the rightmost infected site by time $t$, equals $\Lambda_{1}[0,t)$. From these and noting that $N_{t} \leq R_{t}^{0}+ \bar{x}_{t}+ D_{t}$, the proof is complete by elementary Poisson processes results.
\end{proof}

%\textsc{Comment.} A blaring ommision from the literature are the survival properties of the process . We are neither able to reason in favor of the strengthen of  \ref{} in an equality nor in contradiction. 

%\begin{remark}
%\textup{ Note that the condition $\mu \geq \lambda$ is needed for assuring that monotonicity of the process (in the sense of Proposition \ref{moninit}) and the right endpoints coupling result (namely, Lemma 2.2 in \cite{T}) hold. These are necessary only in the proof of the existence of ?? in \cite{T} and in that of Corrolary \ref{asconv}, but not in the proof of this theorem above.
%The punchline can be proved without this assumption?}
%\end{remark} 

%\textbf{Acknowledgements:} The author is appreciative to Stan Zachary for research advising and helpful discussions, and to Pablo Ferrari for hospitality. %while research for this work was performed. %which was financially supported by a Heriot-Watt University scholarship. 
%Thus  hence 

\end{document}